\documentclass[12pt,fleqn]{article}
\usepackage{amsfonts}
\usepackage{latexsym,amsmath,amssymb}
\date{}

\begin{document}
\title{An improvement upon unmixed decomposition of an algebraic
 variety
\thanks{Partially supported by a NKBRPC (2004CB318000)
}
\author{Zhenyi Ji\thanks{E-mail: jizhenyi0010@163.com (Z.Y.Ji)},
Yongbin Li\thanks{E-mail: yongbinli@uestc.edu.cn (Y.B.Li)}.\\
\emph{\small School of Applied Mathematics,}\\
   \emph{\small University of Electronic Science and Technology of China,}\\
   \emph{\small Chengdu, Sichuan, 610054, China}}
\date{}}
\maketitle

\begin{abstract}
Decomposing an algebraic variety into irreducible or equidimensional
components is a fundamental task in classical algebraic geometry and
has various applications in modern geometry engineering. Several
researchers studied the problem and developed efficient algorithms
using $Gr$\"{o}$bner$ basis method. In this paper, we try to modify
the computation of unmixed decomposition of an algebraic variety
based on improving the computation of $Zero(sat(\mathbb{T}))$, where
$\mathbb{T}$ is a triangular set in $\textbf{K[X]}$.


\ \ \ \ \ \ \ \ \ \ \ \ \ \ \ $Wu's$ characteristic set,U-set.

\end{abstract}

\section{Introduction}
\ \ \ \ \ \ Let \textbf{K} be a field of characteristic 0 and
$\textbf{K}[x_1,x_2,\ldots,x_n]$ (or $\textbf{K[X]}$ for short) the
ring of polynomials in the variables $(x_1,x_2,\ldots,x_n)$ with
coefficients in \textbf{K}. A polynomial set is a finite set
$\mathbb{P}$ of nonzero polynomials in \textbf{K[X]}. The ideal of
$\textbf{K[X]}$ generated by all elements of $\mathbb{P}$ is denoted
by $Ideal(\mathbb{P})$ and the algebraic variety of $\mathbb{P}$ is
denoted by $Zero(\mathbb{P})$. The method of $Gr$\"{o}$bner$ bases
introduced by $Buchberge$$^{[1,2]}$provides a powerful device for
computing a basis of $Ideal(\mathbb{P})$. It is well known that
$Wu$$^{[22]}$ provided an efficient method for constructing a $Wu's$
characteristic set of every polynomial set to compute the variety of
$\mathbb{P}$ in 1978. Therefrom, various algorithms for computing
triangular decomposition of polynomial sets and  $systems$ are
developed by some researchers$^{[3,4,10,11,17,19,23,24,25]}$.

Based on various triangular decompositions for polynomial systems,
including the famous $Wu's$ characteristic set method, and the
$Gr$\"{o}$bner$ Basis method we can get the unmixed decomposition of
an algebraic variety.

According to the analytic method established by $Zhang$
et.al$^{[28]}$, we get the modified method to compute characteristic
series. Furthermore, we try to improve the computation of
$Zero(sat(\mathbb{T}))$, where $\mathbb{T}$ is a triangular set.
Some examples can illustrate our improvement.

\section{Preliminaries}
\subsection{$Wu's$ characteristic set}
\ \ \ \ \ \ For any polynomial $ p \notin \textbf{K}$, the biggest
index $k$ such that $deg(p,x_k)>0$ is called the \emph{class}, $x_k$
the \emph{leading variable}, $deg(p,x_k)$ the \emph{leading degree
}of $p$, and $lcoeff(p,x_k)$ the \emph{leading coefficient} of $p$,
denoted by $cls(p)$, $lv(p)$, $ldeg(p)$, $ini(p)$, respectively.\\

\textbf{Definition 2.1.1.} A finite nonempty ordered set
$\mathbb{T}=[f_1,f_2,\ldots,f_s]$ of polynomials in $
\textbf{K[X]}\backslash \textbf{K} $ is called a \emph{triangular
set} if $cls(f_1)<cls(f_2)<\cdots<cls(f_s)$. Triangular set
$\mathbb{T}$ is written as the following form
\begin{equation}
\mathbb{T}=[f_1(u_1,\ldots,u_r,y_1),\ldots,f_s(u_1,\ldots,u_r,y_1,\ldots,y_s)]
\end{equation}
where $(u_1,\ldots,u_r,y_1,\ldots,y_s)$ is a permutation of
$(x_1,\ldots,x_n)$.\\

Let $f\neq 0$ be a polynomial and $g$ any polynomial in
$\textbf{K[X]}$,the $ \emph{pesudo-}\emph{remainde}r$ of $g$ with
respect to $f$ in $lv(f)$ is denoted by $ prem(g,f,lv(f)) $. One can
find a formal definition of $pesudo-remainder$$^{[7,18]}$or two
alternative ones$^{[18,19]}$. For any polynomial $p$ and triangular
set $\mathbb{T}$ $prem(p,\mathbb{T})$ stands for the
\emph{pesudo-remainder}
 of $p$ with respect to $\mathbb{T}$ is defined by
\begin{align}
prem(p,\mathbb{T})=prem(\ldots prem(p,f_s,y_s),\ldots,f_1,y_1).
\end{align}
One can easily deduce the following $\emph{pesudo-remainder
formula}$
\begin{align}
\prod\limits_{i = 1}^s {ini(f_i )} ^{d_i } p = \sum\limits_{i = 1}^s
{q_i f_i }  + prem(p,\mathbb{T}),
\end{align}
where each $d_i$ is a nonnegative integer and $q_i\in \textbf{K[X]}$
for $1\leq i\leq s$.

For any polynomial set $\mathbb{P}\subset \textbf{K[X]}$, we write

\ \ \ \ $prem(\mathbb{P},\mathbb{T})\triangleq
\{prem(p,\mathbb{T})|p\in \mathbb{P}\}$.\\

Given two polynomials $f$, $g$ $\in \textbf{K[X]}$, the
$\emph{Sylvester resultant}$ of $f$ and $g$ with respect to some
$x_k$ $(1\leq k\leq n)$ is denoted by $res(f,g,x_k)$. Let $p$ be any
polynomial and $\mathbb{T}=[f_1,f_2,\ldots,f_s]$ a triangular set in
$\textbf{K[X]}$ as (1). The polynomial

\ \ \ \ $res(p,\mathbb{T})\triangleq res(\ldots
res(p,f_s,y_s),\ldots,f_1,y_1)$ \\is called the\emph{ resultant} of
$p$ with respect to $\mathbb{T}$.

Let $\mathbb{T}$ is a triangular set as $(1)$ and $p$ any
polynomial, $p$ is said to be reduced with respect to $\mathbb{T}$
if $deg(p,y_i)<deg(f_i,y_i)$ for all $i$. $\mathbb{T}$ is said to be
\emph{noncontradictory ascending set} if every $f\in
\mathbb{T}\cup(ini(\mathbb{T}))$ is reduced to
$\mathbb{T}\setminus\{f\}$.\\

\textbf{Definition 2.1.2.} A triangular set
$\mathbb{T}=[f_1,f_2,\ldots,f_s]$ is called perfect, if
$Zero(\mathbb{T},ini(\mathbb{T}))\neq \emptyset$
.\\

\textbf{Definition 2.1.3.$^{[22]}$} A noncontradictory ascending set
$\mathbb{T}$ is called a $Wu's$ \emph{characteristic set} of
polynomial set
$\mathbb{P}\subset \textbf{K[X]}$ if\\

\ \ \ \ $\mathbb{T}\subset Ideal(\mathbb{P})$,\ \ \ \
$prem(\mathbb{P},\mathbb{T})=\{0\}$.\\

\textbf{Definition 2.1.4.} A finite set
$\mathbb{T}_1,\mathbb{T}_2,\ldots,\mathbb{T}_s$ is called a
\emph{characteristic series} of polynomial set $\mathbb{P}$ in
\textbf{K[X]} if the following zero decomposition holds
\begin{equation}
Zero(\mathbb{P}) = \bigcup\limits_{i = 1}^s
{Zero(\mathbb{T}_i/ini(\mathbb{T}_i))}
\end{equation}
and $prem(\mathbb{P},\mathbb{T}_i)=\{0\}$ for every $i$.\\

\textbf{Defnition 2.1.5.$^{[10,25]}$} A triangular set
$\mathbb{T}=[f_1,f_2,\ldots,f_s]$ is called a \emph{regular set} if
$res(I,\mathbb{T})\neq 0$ for all $I\in ini(\mathbb{T}).$\\

\textbf{Definition 2.1.6.$^{[4,19]}$} A triangular set
$\mathbb{T}=[f_1,f_2,\ldots,f_s]$ as (1) is called a \emph{normal
set} if
$ini(\mathbb{T})\in \textbf{K[U]}$.\\

\textbf{Definition 2.1.7.} Let $\mathbb{T}$ be a triangular set in
\textbf{K[X]}. The \emph{saturation ideal} of $\mathbb{T}$ \\

$sat(\mathbb{T})\triangleq Ideal(\mathbb{T}):J^\infty=\{g\in
\textbf{K[X]}\mid J^q g\in Ideal(\mathbb{T})$ for some $q>0$$\}$,\\
\\where $J=\prod_{f\in\mathbb{T}}ini(f)$.\\

One can compute a basis of $sat(\mathbb{T})$ by the following Lemma.\\

\textbf{Lemma 2.1.8.$^{[5,8,20]}$} Let
$\mathbb{T}=[f_1,f_2,\ldots,f_s]$ be a triangular set in
\textbf{K[X]}, $z$ is a new variable,
$\mathbb{H}=\mathbb{T}\cup\{zJ-1\}=\{f_1,f_2,\ldots,f_s,zJ-1\}$,
$Gb$ be the $Gr$\"{o}$bner$ basis of $\mathbb{H}$ with respect to a
\emph{lexicographic ordering} where $z$ is greater than every $x_i$.
Then
\begin{equation}
sat(\mathbb{T})=Ideal(\mathbb{H})\cap \textbf{K[X]}=Ideal(Gb\cap
\textbf{K[X]}).
\end{equation}
\subsection{The theory of weakly non-degenerate conditions and its application}

\ \ \ \ Let $\mathbb{T}$ be as $(1)$,  we denote ${\mathbb{C}_f}_i$
the set of all the nonzero coefficients of $f_i$ in $y_i$,
$\mathbb{R}_{f_{i}}=\{res(c,\mathbb{T})\neq0:c\in\mathbb{C}_{f_i}\}$
for any $f_i \in\mathbb{T}$. For any
$\overline{\textbf{z}}=({\overline{\textbf{u}}},\overline{y}_1,\ldots,\overline{y}_s)\in
Zero(\mathbb{T})$, we write ${\overline{\textbf{z}}}^{\{{j}\}}$ for
${\overline{\textbf{u}}}$,$\overline{y}_1$,$\ldots$,$\overline{y}_j$
or $(\overline{\textbf{u}},\overline{y}_1,\ldots,\overline{y}_j)$ with $0\leq j\leq s.$ \\

\textbf{Definiton 2.2.1.}$^{[28]}$ Let $\mathbb{T}$ as $(1)$ be a
regular set in $\textbf{K[X]}$. A zero $\textbf{z}_0\in
Zero(\mathbb{T})$ is called a $\emph{quasi-normal zero}$ if $
\textbf{z}_0 ^{\{ i - 1\} } \notin Zero(\mathbb{C}_{f_i } )$ for any
$1\leq i\leq s,$ also said to satisfying the \emph{weakly
non-degenerate conditions}.\\

The following definition is an extension of the
concept of \emph{quasi-normal zero} of regular set to triangular set.\\

\textbf{Definiton 2.2.2.}$^{[13]}$ Let
$\mathbb{T}=[f_1,f_2,\ldots,f_s]$ as (1) be a triangular set in
$\textbf{K[X]}$. A zero \textbf{$z_0$} $\in Zero(\mathbb{T})$ is
called a \emph{quasi-normal zero} of $\mathbb{T}$ if for any $1\leq
i\leq s$, \emph{either conditions holds}:

\ \ \ \ \ \ \ \ a. $I_i({\bf z}^{\{i-1\}}_0)\ne 0 $ if ${
{res}}(I_i,\mathbb{T})=0$;

\ \ \ \ \ \ \ \ b. ${\bf z}^{\{i-1\}}_0 \notin {\rm Zero}( {\mathbb
C}_{f_i})$ if
   ${ res}(I_i,{\mathbb T})\ne 0 $.\\

For any triangular set $\mathbb{T}$ in $ \textbf{K[X]}$, we denote
$QnZero(\mathbb{T})$ the set of all \emph{quasi-normal zeros} of
$\mathbb{T}$ and $\overline{QnZero(\mathbb{T})}^E$ the closure of
$QnZero(\mathbb{T})$ in topological space $\textbf{K}^n$,  where
$\textbf{K}^n$ is induced by follow metric

$|{\bf z}-{\bf z}^*|={\rm max}\{|x_1-x^*_1|,|x_2-x^*_2|,\ldots,
 |x_n-x^*_n|\}$  for any ${\bf z},{\bf z}^*\in \tilde{{\bf
 K}}^n$,
\\ then we have
the following theorem.\\

\textbf{Theorem 2.2.3.$^{[13]}$} For any triangular set
$\mathbb{T}=[f_1(\textbf{u},y_1),\ldots,f_s(\textbf{u},y_1,\ldots,y_s)]$,
we have\\

\ \ \ \ \ \ \ \ \
$Zero(\mathbb{T}/ini(\mathbb{T}))\subseteq\overline{QnZero(\mathbb{T})}^E\subseteq
Zero(sat(\mathbb{T}))$.\\

The following definition plays an important role in this paper.\\

\textbf{Definition 2.2.4.$^{[9,13]}$} Let $\mathbb{T}$ be a
triangular set in $\textbf{K[X]}$, We establish

$\mathbb{U}_\mathbb{T}\triangleq
\{c\in\mathbb{C}_f:res(ini(f),\mathbb{T})\neq0$ and $
\mathbb{R}_f\cap K=\emptyset, f\in\mathbb{T}\}\cup$

\ \ \ \ \ \ \ \ $\{c:res(c,\mathbb{T})=0,c\in ini(\mathbb{T})\}.$\\

\textbf{Remark}: This definition has slightly difference with the
notion in$^{[9,13]}$.

\textbf{Example 2.2.5.} Let $\mathbb{T}=[f_1,f_2,f_3]$, under
$x_1\prec x_2\prec x_3\prec x_4$, where

\ \ \ \ \ \ \ \ \ $f_1=x_1x_2^2+x_2+2x_1^2,$

\ \ \ \ \ \ \ \ \ $f_2=x_2x_3+x_1x_2^2+x_2x_1+2x_1,$

\ \ \ \ \ \ \ \ \ $f_3=x_2x_4^2-x_4-x_2-x_3.$

By above notation, we know

\ \ \ \ \ \ \ \ \ $\mathbb{C}_{f_1}=\{x_1,1,2x_1^2\}$,

\ \ \ \ \ \ \ \ \ $\mathbb{C}_{f_2}=\{x_2,x_1x_2^2+x_1x_2+2x_1\}$,

\ \ \ \ \ \ \ \ \ $\mathbb{C}_{f_3}=\{x_3,-1,-x_2-x_3\}$.

It is easy to see that

\ \ \ \ \ \ \ \ \
$\mathbb{R}_{f_1}=\mathbb{C}_{f_1}=\{x_1,1,2x_1^2\}$,

\ \ \ \ \ \ \ \ \ $\mathbb{R}_{f_2}=\{2x_1^2,{x_{{1}}}^{2} \left(
4\,{x_{{1}}}^{4}-6\,{x_{{1}}}^{3}+2\,{x_{{1}}}^{2 }-2\,x_{{1}}+2
\right) \}$,

\ \ \ \ \ \ \ \ \ $\mathbb{R}_{f_3}=\{2x_1
,-1,4\,{x_{{1}}}^{6}-14\,{x_{{1}}}^{5}+14\,{x_{{1}}}^{4}+2\,{x_{{1}}}^{2}-
2\,x_{{1}} \}$.

Thus $\mathbb{U}_{\mathbb{T}}=\{x_2\}$.

Similarly, one can compute that
$\mathbb{U}_{\mathbb{T}^{*}}=\emptyset$ where
$\mathbb{T}^{*}=[g_1,g_2,g_3]$

\ \ \ \ \ \ \ \ \
$g_1=-x_{{2}}{x_{{3}}}^{2}-x_{{3}}+x_{{1}}{x_{{2}}}^{2}-x_{{2}}x_{{1}}
, $

\ \ \ \ \ \ \ \ \
$g_2=x_{{1}}{x_{{4}}}^{2}+x_{{3}}{x_{{4}}}^{2}+x_{{4}}+x_{{3}}-2\,x_{{2}}+x
_{{2}}x_{{1}}+x_{{1}}{x_{{2}}}^{2}
 , $

\ \ \ \ \ \ \ \ \ $g_3=
x_{{3}}{x_{{5}}}^{2}+2\,x_{{5}}+2\,x_{{1}}{x_{{2}}}^{2}+x_{{2}}x_{{1}}
+2\,x_{{3}} , $

under $x_1\prec x_2\prec x_3\prec x_4\prec x_5$.\\

\textbf{Theorem 2.2.6.$^{[13]}$} For any triangular set
$\mathbb{T}$, we have

\ \ \ \ \ \ \ \ \ $Zero(\mathbb{T}/\mathbb{U}_\mathbb{T})\subseteq
\overline{QnZero(\mathbb{T})}^E\subseteq Zero(sat(\mathbb{T}))$.\\

\textbf{Corollary 2.2.7.$^{[13]}$} Let $\mathbb{T}$ be a triangular
set in $\textbf{K[X]}$ with $\mathbb{U}_\mathbb{T}=\emptyset$. Then

\ \ \ \ \ \ \ \ \  $Zero(\mathbb{T})=Zero(sat(\mathbb{T}))$.\\

\ \ \ \ \ \ Based on the the theory of weakly non-degenerate
conditions we get the following modified algorithm
$CharserA$$^{[9,13]}$ to compute
characteristic series.\\

\textbf{Algorithm CharserA:} $\Psi\leftarrow CharserA(\mathbb{P})$.
Given a nonempty polynomial set $\mathbb{P}$ in \textbf{K[X]}, this
algorithm computes a finite sets $\Psi$ such that

\ \ \ \ \ \ \ \ \ \ \ \ $ Zero(\mathbb{P}) =
\bigcup\limits_{\mathbb{T} \in \Psi }
{Zero(\mathbb{T}/\mathbb{U}_\mathbb{T} )}$

$C1:$ Set $\Phi  \leftarrow \{ \mathbb{P}\}$, $\Psi\leftarrow
\emptyset$.

$C2:$ While $\Phi \neq \emptyset$ do:

\ \ \ $C2.1$. Let $\mathbb{F}$ be an element of $\Psi$ and set
$\Psi\leftarrow \Psi \backslash {\mathbb{F}}$.

\ \ \ $C2.2$ Compute $\mathbb{T}\leftarrow Charset(\mathbb{F})$.

\ \ \ $C2.3$ If $\mathbb{T}$ is noncontradictory, then compute
$\mathbb{U}_{\mathbb{T}}$.

\ \ \ $C2.4$ If $\mathbb{U}_{\mathbb{T}}=\emptyset$, then set $\Psi
\leftarrow \Psi\cup{\mathbb{T}}$.

\ \ \ $C2.5$ If $\mathbb{U}_{\mathbb{T}}\neq\emptyset$, then set

\ \ \ \ \ \ \ \ \ \ \ \ \ $\Psi \leftarrow \Psi \cup {\mathbb{T}}$,
$\Phi \leftarrow \Phi \{\mathbb{F}\cup \mathbb{T}\cup \{I\}: I\in
\mathbb{U}_{\mathbb{T}}\}$.\\

\textbf{Example 2.2.8.} Let $\mathbb{P}=\{p_1,p_2,p_3\}$ in
$\textbf{K}[x_1,x_2,x_3,x_4]$, where

\ \ \ \ \ \ \ \ \ \ \ \ \
$p_1={x_{{3}}}^{2}+2\,x_{{3}}x_{{2}}+x_{{1}}x_{{2}}+2$,

\ \ \ \ \ \ \ \ \ \ \ \ \
$p_2={x_{{1}}}^{3}+2-4\,x_{{3}}{x_{{2}}}^{2}-2\,x_{{2}}{x_{{3}}}^{2}+2\,{x_
{{1}}}^{2}-2\,x_{{2}}$,

\ \ \ \ \ \ \ \ \ \ \ \ \
$p_3=x_{{2}}x_{{3}}{x_{{4}}}^{2}+x_{{4}}+{x_{{1}}}^{2}+2\,x_{{3}}+x_{{2}}$.

Under the variable ordering $x_1\prec x_2\prec x_3 \prec x_4$. By
the above description, one can easily get $CharserA=\{\mathbb{T}\}$,
where

$\mathbb{T}=[2\,x_{{1}}{x_{{2}}}^{2}+2\,x_{{2}}+{x_{{1}}}^{3}+2\,{x_{{1}}}^{2}+2,{x_{{3}}}^{2}+2\,x_{{2}}x_{{3}}+x_{{1}}x_{{2}}+2,x_{{2}}x_{{3}}{x_{{4}}}^{2}+x_{{4}}+$

\ \ \ \ \ \ \ ${x_{{1}}}^{2}+2\,x_{{3}}+x_{{2}}]$.

It is easy to see that $\mathbb{U}_{T}=\emptyset$, this implies that

$Zero(\mathbb{P})=Zero(\mathbb{T})$

Compared with $MMP$, $epsilon$, $Regular$, we get the following zero
decomposition directly.

$Zero(\mathbb{P})= Zero(\mathbb{T}_1/\{x_1,x_2x_3\})
\bigcup\limits_{i = 2}^5 {Zero( \mathbb{T}_i)} $. ($MMP$)

where

$\mathbb{T}_1=[2\,x_{{1}}{x_{{2}}}^{
2}+2\,x_{{2}}+2+2\,{x_{{1}}}^{2}+{x_{{1}}}^{3},
{x_{{3}}}^{2}+2\,x_{{2}}x_{{3}}+x_{{1}}x_{{2}}+2,x_{{2}}x_{{3}}{x_{{4}}}^{2}+x_{{4}}+$

\ \ \ \ \ \ \ \ ${x_{{1}}}^{2}+2\,x_{{3}}+x_{{2}}] .$

$\mathbb{T}_2=[-{x_{{1}}}^{7}-4\,{x_{{1}}}^{6}-4\,{x_{{1}}}^{5}-4\,{x_{{1}}}^{4}-12
\,{x_{{1}}}^{3}-8\,{x_{{1}}}^{2}-4\,x_{{1}}-8,-2\,x_{{2}}+{x_{{1}}}^{3
}+2\,{x_{{1}}}^{2}$

\ \ \ \ \ \ \ \ \ \ \ \
$+2,-{x_{{1}}}^{3}x_{{3}}-2\,{x_{{1}}}^{2}x_{{3}}-2\,
x_{{3}},{x_{{1}}}^{4}x_{{4}}+2\,{x_{{1}}}^{3}x_{{4}}+2\,x_{{1}}x_{{4}}
+{x_{{1}}}^{6}+2\,{x_{{1}}}^{5}-4\,{x_{{1}}}^{2}-4] .$

$\mathbb{T}_3=[-2\,{x_{{1}}}^{3}-4\,{x_{{1}}}^{2}-4,4\,x_{{2}},4\,{x_{{3}}}^{2}+8,4
\,x_{{4}}+8\,x_{{3}}+4\,{x_{{1}}}^{2}] .$

$\mathbb{T}_4=[-{x_{{1}}}^{3}-2\,{x_{{1}}}^{2}-2,-2\,x_{{2}},-2\,{x_{{3}}}^{2}-4,-2
\,x_{{4}}-4\,x_{{3}}-2\,{x_{{1}}}^{2}] .$

$\mathbb{T}_5=[2\,x_{{1}},4\,x_{{2}}+4,2\,{x_{{3}}}^{2}-4\,x_{{3}}+4,-4\,x_{{3}}{x_{
{4}}}^{2}+4\,x_{{4}}+8\,x_{{3}}-4] $.

$Zero(\mathbb{P})=  Zero(\mathbb{T}_1/\{x_1,x_2x_3\})
\bigcup\limits_{i = 2}^4 {Zero( \mathbb{T}_i)}$. ($epsilon$)

where

$\mathbb{T}_1=[2\,x_{{1}}{x_{{2}}}^{2}+2\,x_{{2}}+2+2\,{x_{{1}}}^{2}+{x_{{1}}}^{3},{
x_{{3}}}^{2}+2\,x_{{2}}x_{{3}}+x_{{1}}x_{{2}}+2,x_{{2}}x_{{3}}{x_{{4}}
}^{2}+x_{{4}}$

\ \ \ \ \ \ \ \ $+{x_{{1}}}^{2}+2\,x_{{3}}+x_{{2}}].$

$\mathbb{T}_2=[{x_{{1}}}^{4}+2\,{x_{{1}}}^{3}+2\,x_{{1}}+4,2\,x_{{2}}-{x_{{1}}}^{3}-
2-2\,{x_{{1}}}^{2},x_{{3}},2\,x_{{4}}+4\,{x_{{1}}}^{2}+{x_{{1}}}^{3}+2
] .$

$\mathbb{T}_3=[x_{{1}},1+x_{{2}},{x_{{3}}}^{2}+2-2\,x_{{3}},x_{{3}}{x_{{4}}}^{2}-x_{
{4}}-2\,x_{{3}}+1] . $

$\mathbb{T}_4=[{x_{{1}}}^{3}+2+2\,{x_{{1}}}^{2},x_{{2}},{x_{{3}}}^{2}+2,x_{{4}}+{x_{
{1}}}^{2}+2\,x_{{3}}].$

$Zero(\mathbb{T})=Zero(\mathbb{T}/ \{x_1,x_2x_3\})$. $(Regular)$

where

$\mathbb{T}=[2\,x_{{1}}{x_{{2}}}^{2}+2\,x_{{2}}+2+2\,{x_{{1}}}^{2}+{x_{{1}}}^{3},{
x_{{3}}}^{2}+2\,x_{{2}}x_{{3}}+x_{{1}}x_{{2}}+2,x_{{2}}x_{{3}}{x_{{4}}
}^{2}$

\ \ \ \ \ \ \ $+x_{{4}}+{x_{{1}}}^{2}+2\,x_{{3}}+x_{{2}}] ;$

Other experiment comparison see table 1.
\begin{table}[h]
\setlength{\abovecaptionskip}{0pt}
\begin{center}\caption{Number of triangular sets for nine test examples.}
\begin{tabular}{||r||r||r||r|r||r||r||}
\hline Polynomial set  &\itshape $MMP^{[7]}$ &\itshape
$epsilon^{[21]
}$ &\itshape $Regular^{[16]}$ &\itshape $CharserA^{[9,13]}$  \\
\hline &\itshape Branches
  &\itshape Branches  &\itshape Branches
  &\itshape Branches \\
\hline $\mathbb{P}_1$ &1  &1  &2  &1 \\ \hline $\mathbb{P}_2$ &2 &3
 &5 &2 \\ \hline
 $\mathbb{P}_3$ &3  &3  &6  &1 \\ \hline
  $\mathbb{P}_4$ &3  &2  &1  &2 \\ \hline
$\mathbb{P}_5$ &4 &3  &1  &1\\ \hline
     $\mathbb{P}_6$ &6  &?  &2  &4  \\ \hline
$\mathbb{P}_7$ &9  &5  &7  &1 \\ \hline
        $\mathbb{P}_8$ &9  &6  &1  &2 \\ \hline
         $\mathbb{P}_{9}$ &6  &4  &4 &3 \\ \hline
\end{tabular}\\
 \end{center}
\end{table}

\section{Traditional method for unmixed decomposition}

\ \ \ \ An algebraic variety is said to be unmixed or
equidimensional if all irredundant irreducible components have the
same dimension.

Refer to the zero decomposition (4) which provides a representation
of the variety $Zero(\mathbb{P})$ in terms of its subvarieties
determined by $\mathbb{T}_i$. However, each
$Zero(\mathbb{T}_i/ini(\mathbb{T}_i))$ is not necessarily an
algebraic variety, it is a \emph{quasi-algebraic variety}. In what
follows, we shall see how a corresponding variety decomposition may
be obtained by determining, from each $\mathbb{T}_i$, a finite set
of
polynomials.\\

\textbf{Theorem 3.1} Let $\mathbb{P}$ be a polynomial set in
$\textbf{K[X]}$, $\mathbb{T}_1,\mathbb{T}_2,\ldots,\mathbb{T}_s$ is
a characteristic series of $\mathbb{P}$, then we have
\begin{equation}
Zero(\mathbb{P})=\bigcup_{i=1}^s Zero(sat(\mathbb{T}_i)).
\end{equation}\\
The following result provides a useful criterion for removing some
redundant subvarieties in the
decomposition (4) without computing their defining sets.\\

\textbf{Lemma 3.2}$^{[3]}$ Let $\mathbb{P}$ and $\mathbb{T}_i$ be as
in Theorem $3.1$, if $|\mathbb{T}_i|>|\mathbb{P}|$, then we have
\begin{equation}
Zero(sat(\mathbb{T}_i))\subset\mathop  \bigcup \limits_{\mathop {1
\le i \le s}\limits_{i \ne j} } Zero(sat(\mathbb{T}_j)).
\end{equation}

The next theorem can make sure this
decomposition is unmixed.\\

\textbf{Theorem 3.3}$^{[6]}$ Let $\mathbb{T}=[f_1,f_2,\ldots,f_s]$
be a triangular set in $\textbf{K[X]}$, if $\mathbb{T}$ is not
perfect, then $sat(\mathbb{T})=\textbf{K[X]}$. If $\mathbb{T}$ is
perfect, then $Zero(sat(\mathbb{T}))$ is an unmixed variety of dimension $n-|\mathbb{T}|$.\\

For every $i$ let $\mathbb{G}_i$ be the finite basis of
$sat(\mathbb{T}_i)$, which can computed by Lemma 2.1.6. If
$sat(\mathbb{T}_i)=\textbf{K[X]}$ then
$Zero(sat(\mathbb{T}_{i}))=\emptyset$, hence, we can remove it. Thus
a variety decomposition of the following form is obtained:
\begin{equation}
Zero(\mathbb{P})=\bigcup_{i=1}^s Zero(\mathbb{G}_i)
\end{equation}
By theorem 3.3, each $\mathbb{G}_i$ defines an unmixed algebraic
variety.
\section{Improvement}
\ \ \ \ \ \ Based upon the theory of weakly non-degenerate conditions, we have the following result.\\

\textbf{Theorem 4.1} Let $\mathbb{T}=[f_1,f_2,\ldots,f_s]$ as $(1)$
be a triangular set in $\textbf{K[X]}$, then
\begin{equation}
Zero(sat(\mathbb{T}))=Zero(Ideal({\mathbb{T}}):U^{\infty})
\end{equation}
where $U  = \prod\limits_{u \in \mathbb{U}_{\mathbb{T}} } u.$

\textbf{Proof}: It is easy to see that $Zero(sat(\mathbb{T}))\subset
Zero(Ideal({\mathbb{T}}):U^{\infty}).$

Next, to establish containment in the opposite direction, let

\ \ \ \ \ \ \ \ \ \ \ $\textbf{X}\in
Zero(Ideal({\mathbb{T}}):U^{\infty}).$

Equivalently,

\ \ \ \ \ \ \ \ \ \ \ if $fU^m\in Ideal(\mathbb{T})$ for some $m>0$,
then $f(\textbf{X})=0$.

Now, let $f\in Ideal(Zero(\mathbb{T}/U)$, then $fU$ vanishes on
$Zero(\mathbb{T})$. Thus, by the $Nullstellensatz$, $fU\in
\sqrt{Ideal(\mathbb{T})}$, so $(fU)^l\in Ideal(\mathbb{T})$ for some
$l>0$. Hence, $f^l \in Ideal(\mathbb{T}):U^\infty$. we have
$f(\textbf{X})=0$, Thus

\ \ \ \ \ \ \ \ \ \ \ $\textbf{X}\in
Zero(Ideal(Zero(\mathbb{T}/U))$.

This establish that

\ \ \ \ \ \ \ \ \ \ \ $Zero(Ideal(\mathbb{T}):U^\infty)\subset
Zero(Ideal(Zero(\mathbb{T}/U))$.
\\We also know

\ \ \ \ \ \ \ \ \ \ \ $Zero(Ideal(\mathbb{T}):U^\infty)\supset
Zero(Ideal(Zero(\mathbb{T}/U))$,
\\then

$Zero(Ideal(\mathbb{T}):U^\infty)= Zero(Ideal(Zero(\mathbb{T}/U))$

\ \ \ \ \ \ \ \ \ \ \ \ \ \ \ \ \ \ \ \ \ \ \ \ \ \ \ \
$=\overline{Zero(\mathbb{T}/U)}.$
\\where $\overline{Zero(\mathbb{T}/U)}$ denotes the
$\emph{Zariski}$ \emph{closure} of $Zero(\mathbb{T}/U)$.

From the definition of $\emph{Zariski}$ \emph{closure} and
$Zero(\mathbb{T}/U)\subset Zero(sat(\mathbb{T}))$, then
$Zero(sat(T))\supset Zero(Ideal({\mathbb{T}}):U^{\infty}) $.

This completes the proof.

Applying $CharserA^{[9,13]}$ and above theorem we get the
modified algorithm for unmixed decomposition of an algebraic variety.\\

\textbf{Algorithm 4.2}: $\Psi \leftarrow UnmVarDec(\mathbb{P})$.
Given a nonempty set $\mathbb{P}$, this algorithm computes finite
set $\Psi$ of polynomial sets $\mathbb{G}_{1},\ldots,\mathbb{G}_s$
such that the decomposition (8) holds and each $\mathbb{G}_i$
defines an unmixed algebraic variety.

\textbf{U1}: Compute $\Phi\leftarrow CharserA(\mathbb{P})$, and
 set $\Psi\leftarrow \varnothing$:

\textbf{U2}: While $\Phi\neq \emptyset$, do:

 \ \ \ \textbf{U2.1}: Let $\mathbb{T}$ be an element in $\Phi$, and
 $\Phi\leftarrow\Phi\setminus\{\mathbb{T}\}$. if $|\mathbb{T}|>|\mathbb{P}|$, then go to

 \ \ \ \ \ \ \ \ \ \ \ \ $U2$:

 \ \ \ \textbf{U2.2}: Compute $Gr$\"{o}$bner$ basis $\mathbb{G}$ of $Ideal(\mathbb{T}):U^\infty$ according to

  \ \ \ \ \ \ \ \ \ \ \ \ Lemma 2.1.8, and
 set $\Psi\leftarrow\Psi\cup\{\mathbb{G}\}$:

\textbf{U3}: While $\exists $ $\mathbb{G},\mathbb{G}^{*}$ such that
 $rem(\mathbb{G},\mathbb{G}^*)=\{0\}$,do:

  \ \ \ \ \ \ \ \ \ \ \ \ set $\Psi\leftarrow \Psi \
  \backslash \{\mathbb{G}^*\}$.\\
\\ where
$rem(\mathbb{G},\mathbb{G}^*)\triangleq\{rem(p,\mathbb{G}^*)|p\in
\mathbb{G}\}$ and $rem(p,\mathbb{G}^*)$ see$^{[1,2]}$ for details.\\

\textbf{Example 4.3.} Let $\mathbb{P}=\{f_1,f_2,f_3,f_4\}$ be a
polynomial set in $\textbf{K}[x_1,x_2,x_3,x_4]$ where

$g_1 =
2\,{x_{{3}}}^{2}x_{{1}}+{x_{{3}}}^{2}x_{{2}}+x_{{3}}+x_{{1}},$

$g_2=-x_{{2}}x_{{3}}x_{{4}}+2\,x_{{1}}{x_{{2}}}^{2}-x_{{3}}{x_{{4}}}^{2}+x_
{{2}}+2\,x_{{1}}-x_{{4}} ,$

$g_3={x_{{5}}}^{2}{x_{{1}}}^{2}-x_{{2}}{x_{{5}}}^{2}+x_{{5}}-x_{{3}}x_{{4}}
+x_{{2}}+x_{{1}}
 ,$

$g_4=2\,{x_{{2}}}^{3}{x_{{3}}}^{2}+2\,x_{{2}}{x_{{3}}}^{3}x_{{4}}+2\,{x_{{4
}}}^{2}{x_{{3}}}^{3}+2\,{x_{{2}}}^{2}x_{{3}}+x_{{2}}x_{{3}}x_{{4}}+2\,
x_{{4}}{x_{{3}}}^{2}+x_{{3}}{x_{{4}}}^{2}-x_{{2}}+$

\ \ \ \ \ \ $2\,x_{{3}}+x_{{4}} .$

Under variable ordering $x_1\prec x_2\prec x_3 \prec x_4\prec x_5$
$\mathbb{P}$ is decomposed into six characteristic sets
$\mathbb{T}_i^*$ such that \ \ \ \ \ \ \ \ \ \ \ \ \ \[
Zero(\mathbb{P}) = \bigcup\limits_{i = 1}^6
{Zero(\mathbb{T}_i^*/ini(\mathbb{T}_i^* ))}
\] where

$\mathbb{T}^*_1=[2\,x_{{1}}{x_{{3}}}^{2}+x_{{2}}{x_{{3}}}^{2}+x_{{3}}+x_{{1}},-x_{{3}}
{x_{{4}}}^{2}-x_{{2}}x_{{3}}x_{{4}}-x_{{4}}+x_{{2}}+2\,x_{{1}}+2\,x_{{
1}}{x_{{2}}}^{2},-{x_{{1}}}^{2}{x_{{5}}}^{2}
 $

\ \ \ \ \ \ \ \ $+x_{{2}}{x_{{5}}}^{2}-x_{{
5}}+x_{{3}}x_{{4}}-x_{{2}}-x_{{1}}],$

$\mathbb{T}^*_2=[x_{{2}}-{x_{{1}}}^{2},2\,x_{{1}}{x_{{3}}}^{2}+{x_{{1}}}^{2}{x_{{3}}}^
{2}+x_{{3}}+x_{{1}},x_{{3}}{x_{{4}}}^{2}+x_{{4}}+x_{{3}}{x_{{1}}}^{2}x
_{{4}}-2\,x_{{1}}-{x_{{1}}}^{2}- $

\ \ \ \ \ \ \ \ $2\,{x_{{1}}}^{5},-x_{{5}}+x_{{3}}x_{{4
}}-x_{{1}}-{x_{{1}}}^{2}],$

$\mathbb{T}^*_3=[x_{{2}}+2\,x_{{1}},x_{{3}}+x_{{1}},x_{{1}}{x_{{4}}}^{2}-x_{{4}}-2\,{x
_{{1}}}^{2}x_{{4}}+8\,{x_{{1}}}^{3},{x_{{1}}}^{2}{x_{{5}}}^{2}+2\,x_{{
1}}{x_{{5}}}^{2}+x_{{5}}- $

\ \ \ \ \ \ \ \ $x_{{1}}+x_{{1}}x_{{4}}],$

$\mathbb{T}^*_4=[x_{{1}},x_{{3}},-x_{{4}}+x_{{2}},x_{{2}}{x_{{5}}}^{2}-x_{{5}}-x_{{2}}
]
 ,$

$\mathbb{T}^*_5=[x_{{1}}+2,x_{{2}}-4,x_{{3}}-2,2\,{x_{{4}}}^{2}+9\,x_{{4}}+64,x_{{5}}+
2-2\,x_{{4}}]
 ,$

$\mathbb{T}^*_6=[x_{{1}},x_{{2}},x_{{3}},x_{{4}},x_{{5}}].$

 $\mathbb{T}^*_5$ and $\mathbb{T}^*_6$ contains five polynomials and thus need not be
considered for the variety decomposition by Lemma 3.2. In order to
obtain an unmixed decomposition of $Zero(\mathbb{P})$, It remains to
determine $sat(\mathbb{T}^*_1)$, $sat(\mathbb{T}^*_2)$,
$sat(\mathbb{T}^*_3)$, $sat(\mathbb{T}^*_4)$ by computing the
respectively $Gr$\"{o}$bner$ basis $\mathbb{G}_1$, $\mathbb{G}_2$,
$\mathbb{G}_3$, $\mathbb{G}_4$ of $\mathbb{T}^*_1\cup
\{1+zx_3(2x_1+x_2)(x_2-x_1^2)\}$,
$\mathbb{T}_2^*\cup\{1-zx_3(2x_1+x_1^2)\}$,
$\mathbb{T}_3^*\cup\{1-zx_1(x_1+2)\}$,
$\mathbb{T}_4^*\cup\{1-zx_2\}$ according to Lemma 2.1.8. The
$Gr$\"{o}$bner$ base may be found to consist of 15, 9, 10 and 6
polynomials respectively. Let $sat(\mathbb{T}^*_i)=\mathbb{G}_i\cap
\textbf{K}[x_1,\ldots,x_5]$ for $i=1, 2, 3, 4$. We have

$sat(\mathbb{T}^*_1)= $

$\left[ {\begin{array}{*{20}c}
{2\,{x_{{3}}}^{2}x_{{1}}+x_{{2}}{x_{{3}}}^{2}+x_{{3}}+x_{{1}},}
\hfill \\
{2\,x_{{3}}{x_{{2}}}^{3}x_{{1}}+4\,x_{{3}}{x_{{2}}}^{2}{x_{{1}}}^{2}-x_
{{4}}x_{{3}}x_{{2}}+x_{{3}}{x_{{2}}}^{2}+x_{{1}}{x_{{4}}}^{2}-2\,x_{{1
}}x_{{4}}x_{{3}}+x_{{1}}x_{{4}}x_{{2}}+}\hfill \\
{4\,x_{{1}}x_{{3}}x_{{2}}+2\,{x_
{{2}}}^{2}x_{{1}}+4\,x_{{3}}{x_{{1}}}^{2}-x_{{4}}+x_{{2}}+2\,x_{{1}}
,} \hfill \\
{x_{{4}}x_{{3}}x_{{2}}-2\,{x_{{2}}}^{2}x_{{1}}+{x_{{4}}}^{2}x_{{3}}-x_{
{2}}-2\,x_{{1}}+x_{{4}} ,}\hfill \\
{-{x_{{5}}}^{2}{x_{{1}}}^{2}+x_{{2}}{x_{{5}}}^{2}-x_{{5}}+x_{{4}}x_{{3}
}-x_{{2}}-x_{{1}},
}\hfill \\
{{x_{{5}}}^{2}x_{{3}}-{x_{{3}}}^{2}x_{{1}}-x_{{3}}-x_{{1}}+{x_{{5}}}^{2
}{x_{{3}}}^{2}{x_{{1}}}^{2}+x_{{5}}{x_{{3}}}^{2}+x_{{1}}{x_{{5}}}^{2}-
{x_{{3}}}^{3}x_{{4}}+ }\hfill
\\
{2\,{x_{{5}}}^{2}{x_{{3}}}^{2}x_{{1}}.}\hfill \\
\end{array}} \right]
$

$sat(\mathbb{T}_2^*)=$

$\left[ {\begin{array}{*{20}c} {-x_1+x_2^2,}\hfill\\
{2\,{x_{{3}}}^{2}x_{{1}}+x_{{3}}+x_{{1}}+{x_{{1}}}^{2}{x_{{3}}}^{2},}\hfill
\\
{2\,x_{{3}}{x_{{2}}}^{3}x_{{1}}+4\,x_{{3}}{x_{{2}}}^{2}{x_{{1}}}^{2}-x_
{{3}}x_{{2}}x_{{4}}+x_{{3}}{x_{{2}}}^{2}+{x_{{4}}}^{2}x_{{1}}-2\,x_{{4
}}x_{{3}}x_{{1}}+x_{{4}}x_{{2}}x_{{1}}+
}\hfill\\
{4\,x_{{3}}x_{{2}}x_{{1}}+2\,x_{
{1}}{x_{{2}}}^{2}+4\,x_{{3}}{x_{{1}}}^{2}-x_{{4}}+x_{{2}}+2\,x_{{1}},}\hfill
\\
{x_{{3}}x_{{2}}x_{{4}}-2\,x_{{1}}{x_{{2}}}^{2}+x_{{3}}{x_{{4}}}^{2}-x_{
{2}}-2\,x_{{1}}+x_{{4}}, }\hfill \\
{-{x_{{5}}}^{2}{x_{{1}}}^{2}+x_{{2}}{x_{{5}}}^{2}-x_{{5}}+x_{{3}}x_{{4}
}-x_{{2}}-x_{{1}} ,}\hfill  \\
{-{x_{{3}}}^{2}x_{{1}}-x_{{3}}-x_{{1}}+{x_{{5}}}^{2}x_{{1}}+{x_{{1}}}^{
2}{x_{{5}}}^{2}{x_{{3}}}^{2}+{x_{{3}}}^{2}x_{{5}}-{x_{{3}}}^{3}x_{{4}}
+x_{{3}}{x_{{5}}}^{2}}\hfill \\
{+2\,{x_{{5}}}^{2}{x_{{3}}}^{2}x_{{1}} .}\hfill  \\

\end{array}} \right]$

$sat(\mathbb{T}_3^*)=$

$\left[ {\begin{array}{*{20}c} {2x_1+x_2,}\hfill \\
{x_3+x_1,}\hfill \\
{2\,{x_{{3}}}^{2}x_{{1}}+x_{{3}}+x_{{1}}+{x_{{1}}}^{2}{x_{{3}}}^{2},}\hfill
\\
{{x_{{5}}}^{2}{x_{{1}}}^{2}-x_{{1}}+x_{{5}}+2\,{x_{{5}}}^{2}x_{{1}}+x_{
{4}}x_{{1}} ,}\hfill  \\
{{x_{{5}}}^{2}x_{{4}}x_{{1}}+2\,{x_{{5}}}^{2}x_{{4}}+x_{{5}}{x_{{4}}}^{
2}-2\,x_{{5}}x_{{4}}x_{{1}}+8\,x_{{5}}{x_{{1}}}^{2}+{x_{{4}}}^{2}-x_{{
4}} ,}\hfill  \\
{2\,{x_{{4}}}^{2}{x_{{5}}}^{2}+{x_{{4}}}^{3}x_{{5}}+4\,x_{{5}}x_{{4}}{x
_{{1}}}^{2}+16\,x_{{5}}{x_{{1}}}^{3}+9\,{x_{{5}}}^{2}x_{{4}}+4\,x_{{5}
}{x_{{4}}}^{2}+32\,x_{{5}}{x_{{1}}}^{2}-}\hfill
\\
{32\,{x_{{5}}}^{2}x_{{1}}-8\,x_{{5}}x_{{4}
}x_{{1}}+{x_{{4}}}^{3}+4\,x_{{4}}{x_{{1}}}^{2}+16\,{x_{{1}}
}^{3}-4\,x_{{5}}x_{{4}}+8\,x_{{5}}x_{{1}}-6\,x_{{4}}}\hfill
\\
{+3\,{x_{{4}}}^{2}-14\,x_{{4}}
x_{{1}}-8\,{x_{{1}}}^{2}-16\,x_{{5}}+16\,x_{{1}} .}\hfill
\\

\end{array}} \right]$

$sat(\mathbb{T}_4^*)=\mathbb{T}_4^*.$

It is easy to verify that $Zero(sat(\mathbb{T}_2^*))$,
$Zero(sat(\mathbb{T}_3^*))$ and $Zero(sat(\mathbb{T}_4^*))$ are
subvarieties of $Zero(sat(\mathbb{T}_1^*))$. Therefore,
$Zero(\mathbb{P})=Zero(sat(\mathbb{T}^*_1))$ is an unmixed
decomposition.

By our improvement, $\mathbb{P}$ is decomposition into two
characteristic sets such that
$Zero(\mathbb{P})=Zero(\mathbb{T}_1^*/{x_3})\cup
Zero(\mathbb{T}_4^*)$ since $\mathbb{U}_{\mathbb{T}_1^*}=\{x_3\}$
and $\mathbb{U}_{\mathbb{T}_4^*}=\emptyset$, where $\mathbb{T}_1^*$
and $\mathbb{T}_4^*$ as above. In order to determine
$Zero(sa(\mathbb{T}_1^*))$, we only compute the $Gr$\"{o}$bner$ base
$\mathbb{G}_1$ of $\mathbb{T}_1^*\cup \{1-zx_3\}$ according to
Theorem 4.1. The $Gr$\"{o}$bner$ base may be found to consist of 8
polynomials. Let $Ideal(\mathbb{T}^*_1):x_3^\infty=\mathbb{G}_1\cap
\textbf{K}[x_1,\ldots,x_5]$. We have
$Ideal(\mathbb{T}^*_1):x_3^\infty=sat(\mathbb{T}_1^*)$, and remove
$\mathbb{T}_4^*$ according to $U3$, then we get the result as above.

\textbf{Example 4.4.} Let $\mathbb{P}=\{f_1,f_2,f_3,f_4\}$ be a
polynomial set in $\textbf{K}[x_1,x_2,x_3,x_4,x_5]$ where

$f_1=2\,{x_{{2}}}^{2}x_{{1}}+x_{{1}}x_{{2}}+{x_{{5}}}^{2}x_{{3}}+2\,x_{{5}}
+2\,x_{{3}} ,$

$f_2={x_{{3}}}^{2}x_{{2}}+x_{{3}}+2\,x_{{1}}x_{{2}}+3\,{x_{{4}}}^{2}{x_{{1}
}}^{2}+{x_{{4}}}^{2}x_{{2}}+2\,x_{{4}}
 ,$

$f_3=x_{{3}}{x_{{5}}}^{2}+2\,x_{{5}}+2\,x_{{1}}{x_{{2}}}^{2}+3\,x_{{3}}+{x_
{{3}}}^{2}x_{{2}}+2\,x_{{1}}x_{{2}}+{x_{{2}}}^{3}x_{{1}} , $

$f_4=2\,x_{{5}}+3\,x_{{1}}x_{{2}}+x_{{3}}{x_{{5}}}^{2}+4\,x_{{3}}+2\,x_{{2}
}{x_{{3}}}^{2} .$

 Under variable ordering
$x_1\prec x_2\prec x_3 \prec x_4\prec x_5$ $\mathbb{P}$ is
decomposed into eight characteristic sets $\mathbb{T}_i$ such that \
\ \ \ \ \ \ \ \ \ \ \ \ \[ Zero(\mathbb{P}) = \bigcup\limits_{i =
1}^8 {Zero(\mathbb{T}/ini(\mathbb{T}_i ))}
\]
where

$\mathbb{T}_1=[-x_{{2}}{x_{{3}}}^{2}-x_{{3}}+{x_{{2}}}^{2}x_{{1}}-x_{{1}}x_{{2}},-x_{{1}}{x_{{4}}}^{2}-x_{{3}}{x_{{4}}}^{2}-x_{{4}}-x_{{3}}+2\,x_{{2}}-
x_{{1}}x_{{2}}$

\ \ \ \ \ \ \ \
$-{x_{{2}}}^{2}x_{{1}},x_{{3}}{x_{{5}}}^{2}+2\,x_{{5}}+2\,{x_{{2}}}^{2}x_{{1}}+x_{{1}}x_{{2}}
+2\,x_{{3}} ] ,$

$\mathbb{T}_2=[x_{{1}}{x_{{2}}}^{2}-x_{{1}}x_{{2}},x_{{3}},x_{{1}}{x_{{4}}}^{2}+x_{{4}}-2\,x_{{2}}+2\,x_{{1}}x_{{2}},2\,x_{{5}}+3\,x_{{1}}x_{{2}}]
,$

$\mathbb{T}_3=[x_1,-x_3,-x_4+2x_2,2x_5],$

$\mathbb{T}_4=[x_{{1}}{x_{{2}}}^{2}-x_{{1}}x_{{2}}-{x_{{1}}}^{2}x_{{2}}+x_{{1}},-x_3-x_1,x_{{4}}-2\,x_{{2}}+2\,x_{{1}}x_{{2}}-2\,x_{{1}}+{x_{{1}}}^{2}x_{{2}},-x_{{1}}{x_{{5}}}^{2}+$

\ \ \ \ \ \ \ \ $2\,x_{{5}}+3\,x_{{1}}x_{{2}}-4\,x_{{1}}+2\,{x_{{
1}}}^{2}x_{{2}} ] ,$

$\mathbb{T}_5=[-x_1,-x_3,x_4-2x_2,2x_5],$

$\mathbb{T}_6=[x_1,-x_3,x_4-2x_2,2x_5], $

$\mathbb{T}_7=[-x_2,-x_3,-x_1x_4^2-x_4,2x_5],$

$\mathbb{T}_8=[-x_1,-x_2,-x_3,-x_4,2x_5].$

We can remove $\mathbb{T}_8$ according to Lemma 3.2, and remove
$\mathbb{T}_i(i=2,\ldots,7)$ by compute
$sat(\mathbb{T}_i)(i=2,\ldots,7)$ and $U.3$, where

$sat(\mathbb{T}_1)=$

$\left[ {\begin{array}{*{20}c}
   {{  -{x_{{2}}}^{2}x_{{1}}+x_{{3}}+x_{{2}}{x_{{3}}}^{2}+x_{{1}}x_{{2}}}}, \hfill  \\
   {{  -2\,{x_{{2}}}^{2}x_{{3}}+4\,{x_{{2}}}^{2}x_{{1}}-x_{{1}}x_{{2}}x_{{3}}
          +x_{{3}}x_{{1}}{x_{{2}}}^{2}+x_{{2}}x_{{3}}x_{{4}}-{x_{{1}}}^{2}{x_{{2
          }}}^{2}-{x_{{1}}}^{2}{x_{{2}}}^{3}+}}
          \hfill  \\
   {{  x_{{3}}x_{{1}}{x_{{2}}}^{3}-{x_{{4}}}^{2}{x_{{1}}}^{2}x_{{2}}+{x_{{4}}}^{2}x_{{1}}{x_{{2}}}^{2}-{
          x_{{4}}}^{2}x_{{1}}x_{{2}}-x_{{4}}x_{{1}}x_{{2}}+x_{{1}}{x_{{4}}}^{2}+
          x_{{4}}-2\,x_{{2}}, }} \hfill \\
   {{   x_{{1}}{x_{{4}}}^{2}+{x_{{4}}}^{2}x_{{3}}+x_{{4}}+x_{{3}}-2\,x_{{2}}+x
             _{{1}}x_{{2}}+{x_{{2}}}^{2}x_{{1}}, }} \hfill \\
   {{  2\,x_{{5}}+{x_{{5}}}^{2}{x_{{2}}}^{2}x_{{1}}+2\,x_{{3}}x_{{1}}{x_{{2}}
          }^{3}-{x_{{5}}}^{2}x_{{1}}x_{{2}}+x_{{3}}x_{{1}}{x_{{2}}}^{2}+2\,x_{{2
          }}x_{{3}}x_{{5}}+4\,{x_{{2}}}^{2}x_{{1}}
          }} \hfill  \\
   {{-x_{{1}}x_{{2}}}}, \hfill \\
   {{  2\,{x_{{2}}}^{2}x_{{1}}+x_{{1}}x_{{2}}+{x_{{5}}}^{2}x_{{3}}+2\,x_{{5}}
         +2\,x_{{3}},
          }} \hfill  \\
   {{ -2\,{x_{{4}}}^{2}x_{{5}}-x_{{2}}x_{{4}}+2\,{x_{{2}}}^{2}-8\,{x_{{2}}}^
{3}-2\,{x_{{2}}}^{3}x_{{3}}+4\,x_{{4}}{x_{{2}}}^{2}-4\,{x_{{2}}}^{4}x_
{{3}} -8\,x_{{5}}{x_{{2}}}^{2}
}} \hfill  \\
   {{-2\,x_{{5}}{x_{{2}}}^{3}x_{{3}}-2\,x_{{5}}{x_{{2}}}^{2}x_{{3}}-2\,{x_{
{5}}}^{2}{x_{{2}}}^{3}+2\,{x_{{5}}}^{2}{x_{{2}}}^{2}+{x_{{5}}}^{2}x_{{
4}}{x_{{2}}}^{2}+2\,x_{{4}}{x_{{2}}}^{3}x_{{3}}
}} \hfill  \\
   {{ -{x_{{5}}}^{2}x_{{4}}x_{{2}}+{x_{{2}}}^{2}x_{{4}}x_{{3}}-2\,x_{{5}}{x_
{{2}}}^{2}{x_{{4}}}^{2}+2\,{x_{{4}}}^{2}x_{{1}}x_{{5}}x_{{2}}+2\,x_{{5
}}x_{{4}}x_{{2}}+2\,x_{{2}}x_{{3}}x_{{5}}
}} \hfill  \\
   {{  +2\,x_{{5}}x_{{1}}{x_{{2}}}^{2}+2\,x_{{5}}x_{{1}}{x_{{2}}}^{3}+2\,x_{{5
}}{x_{{4}}}^{2}x_{{2}},
}} \hfill  \\
   {{-8\,x_{{2}}+4\,x_{{4}}-4\,x_{{5}}+2\,x_{{1}}x_{{2}}+4\,x_{{1}}{x_{{4}}
}^{2}+{x_{{5}}}^{2}x_{{4}}-4\,{x_{{2}}}^{2}x_{{3}}-2\,{x_{{4}}}^{2}x_{
{5}}+
}} \hfill  \\
   {{ 4\,{x_{{2}}}^{2}x_{{1}}-2\,x_{{1}}x_{{2}}x_{{3}}+x_{{3}}x_{{1}}{x_{{2}}}^{2}+2\,x_{{2}}x_{{3}
}x_{{4}}-2\,{x_{{1}}}^{2}{x_{{2}}}^{2}-2\,{x_{{1}}}^{2}{x_{{2}}}^{3}+{
x_{{4}}}^{2}{x_{{5}}}^{2}x_{{1}}
}} \hfill  \\
   {{ -2\,x_{{2}}x_{{3}}x_{{5}}-2\,{x_{{4}}}^{2}{x_{{1}}}^{2}x_{{2}}-3\,{x_{{4}}}^{2}x_{{1}}x_{{2}}-2
\,x_{{4}}x_{{1}}x_{{2}}+2\,{x_{{5}}}^{2}x_{{1}}x_{{2}}-2\,{x_{{5}}}^{2
}x_{{2}}
}}. \hfill  \\
\end{array}} \right]
$

then $Zero(\mathbb{P})=Zero(sat(\mathbb{T}_1))$.

By our improvement we get $CharserA(\mathbb{P})=\{\mathbb{T}^{*}\}$,
where

$\mathbb{T}^{*}= \left[ {\begin{array}{*{20}c}
 {-x_{{2}}{x_{{3}}}^{2}-x_{{3}}+x_{{1}}{x_{{2}}}^{2}-x_{{1}}x_{{2}},}\hfill
 \\
 {x_{{1}}{x_{{4}}}^{2}+x_{{3}}{x_{{4}}}^{2}+x_{{4}}+x_{{3}}-2\,x_{{2}}+x
_{{1}}x_{{2}}+x_{{1}}{x_{{2}}}^{2}, }\hfill  \\
{
x_{{3}}{x_{{5}}}^{2}+2\,x_{{5}}+2\,x_{{1}}{x_{{2}}}^{2}+x_{{1}}x_{{2}}
+2\,x_{{3}}. }\hfill  \\
\end{array}} \right]
$

It is easy to see that $\mathbb{U}_{\mathbb{T}^*}=\emptyset$, then
$Zero(sat(\mathbb{T}^*))=Zero(\mathbb{T}^*)$ according to Lemma 4.1,
so we get $Zero(\mathbb{P})=Zero(\mathbb{T}^*)$ directly.
 {\small \small}

\end{document}